\documentclass[12pt,a4paper]{amsart}
\usepackage{amsfonts}
\usepackage{amsthm}
\usepackage{amsmath}
\usepackage{amscd}
\usepackage[latin2]{inputenc}
\usepackage{t1enc}
\usepackage[mathscr]{eucal}
\usepackage{indentfirst}
\usepackage{graphicx}
\usepackage{graphics}
\usepackage{pict2e}
\usepackage{epic}
\numberwithin{equation}{section}
\usepackage[margin=2.9cm]{geometry}
\usepackage{epstopdf} 
\usepackage{amssymb}
\usepackage{alltt}
\usepackage{amsmath,mathtools}

\theoremstyle{plain}
\newtheorem{Th}{Theorem}[section]

\theoremstyle{definition}

\newtheorem{Conj}[Th]{Conjecture}

\newtheorem{?}[Th]{Problem}

\begin{document}
\title{A Note on the Bateman-Horn Conjecture}


\author{Weixiong Li}
\address{3813 Patty Berg Ct, Woodridge, IL 60517}
\curraddr{}
\email{weixiongli@comcast.net}
\thanks{The author would like to thank William A. Stein and the Sage Development Team for making the software freely available to the public. This work would not have been possible without the Sage Mathematics Software. }


\date{}

\dedicatory{}

\begin{abstract}
We report the results of our empirical investigations on the Bateman-Horn conjecture. This conjecture, in its commonly known form, produces rather large deviations when the polynomials involved are not monic. We propose a modified version of the conjecture which empirically demonstrates remarkable accuracy even for modest values of primes.
\end{abstract}

\maketitle

\section{Introduction}
In 1962, Bateman and Horn \cite{BH01}\cite{BH02} proposed the following general conjecture concerning the distribution of primes generated from a set of polynomials\footnote{The wording of the conjecture is mostly from \cite{BH02}.}.

\begin{Conj}[The Bateman-Horn Conjecture]

Let $f_{1},f_{2},\cdots,f_{M}\in\mathbb{Z}[x]$ be distinct irreducible polynomials with positive leading coefficients, and let

\begin{equation}
\pi(f_{1},f_{2},\cdots,f_{M};x)=\#\{n\leq x: f_{1}(n),f_{2}(n),\cdots,f_{M}(n) \text{ are primes}\}
\end{equation}
Suppose that $f(n)=\prod_{i=1}^{M}f_{i}(n)$ does not vanish identically modulo any prime, then for large values of $x$ we have the following asymptotic expression 

\begin{equation}\label{Bate_Horn_Conj}
\pi(f_{1},f_{2},\cdots,f_{M};x)\sim \frac{C(f_{1},f_{2},\cdots,f_{M})}{\prod_{i=1}^{M}\deg f_{i}}\int_{2}^{x}\frac{dt}{(\log t)^{M}}
\end{equation}
in which
\begin{equation}\label{Bateman_Horn_Const}
C(f_{1},f_{2},\cdots,f_{M})=\prod_{p}\Big(1-\frac{1}{p}\Big)^{-M}\Big(1-\frac{\omega_{f}(p)}{p}\Big)
\end{equation}
where the infinite product is over all primes $p$, and $\omega_{f}(p)$ is the number of solutions to the congruence equation $f(n)\equiv 0\, (\text{mod } p)$.
\end{Conj}

The Bateman-Horn conjecture is very general, and many well-known conjectures, such as the Hardy-Littlewood Conjectures B, E, F, K, X, P \cite{HL1}\cite{HL2}, are all special cases of this conjecture.

In proposing this conjecture which now bears their names, Bateman and Horn provided the following heuristic arguments \cite{BH01}. Since the probability of any number $x$ being a prime number is roughly $\frac{1}{\log x}$, the probability of $f_{i}(n)$ being a prime number is $\frac{1}{\log f_{i}(n)}$, which for large $n$ can be approximated by $\frac{1}{(\deg f_{i})\log n}$, hence the expression in equation (\ref{Bate_Horn_Conj}).

Two questions naturally arise. (1) How good is the approximation of replacing $\frac{1}{\log f_{i}(n)}$ with $\frac{1}{(\deg f_{i})\log n}$, and (2) how much better would the 'Bateman-Horn Conjecture' be without this approximation? To the best of this author's knowledge, there seems to be little discussion of these two questions in the literature \cite{BH02}.

Let us perform a naive estimation to answer the first question. Let $f_{i}(n)=a_{i}n^{m_{i}}+b_{i}n^{m_{i}-1}+\cdots$, then $\log f_{i}(n)=\log (a_{i}n^{m_{i}})+\log(1+o(\frac{1}{n}))=m_{i}\log n +\log a_{i}+o(\frac{1}{n})$. we therefore see that when the polynomial is not monic, the relative errors caused by ignoring the term $\log a_{i}$ could be potentially significant for typical counting values around $n\sim 10^{15}$. It is interesting to note that historically most of empirical studies have been performed on monic polynomials where $\log a_{i} = 0$, therefore the errors never showed up in those studies \cite{Shanks1}\cite{Shanks2}. The Hardy-Littlewood Conjecture F was expressed as a function of primes rather than the variable $n$ \cite{HL1}\cite{HL2}, therefore the issue with non-monic polynomials is avoided.

In the next section, we will propose a modified version of the Bateman-Horn conjecture, we then perform empirical computations to compare the two versions of the conjecture with our numerical results.

\section{Empirical Study}

\subsection{The Modified Conjecture}
We will propose the following modified conjecture. 

\begin{Conj}[The Modified Bateman-Horn Conjecture]\label{Modified-Conj}

Let $f_{1},f_{2},\cdots,f_{M}\in\mathbb{Z}[x]$ be distinct irreducible polynomials with positive leading coefficients, and let

\begin{equation}
\pi(f_{1},f_{2},\cdots,f_{M};x)=\#\{n\leq x: f_{1}(n),f_{2}(n),\cdots,f_{M}(n) \text{ are primes}\}
\end{equation}
Suppose that $f(n)=\prod_{i=1}^{M}f_{i}(n)$ does not vanish identically modulo any prime, and $n_{0}$ is the smallest integer such that $\forall\, n>n_{0}$, $f_{i}(n)>1$, $1\leq i \leq M$, then for large values of $x$ we have the following asymptotic expression 

\begin{equation}\label{Modified_Conj_Equation}
\pi(f_{1},f_{2},\cdots,f_{M};x)\sim C(f_{1},f_{2},\cdots,f_{M}) \int_{n_{0}}^{x}\frac{dt}{\prod_{i=1}^{M}\log f_{i}(t)}
\end{equation}
where $C(f_{1},f_{2},\cdots,f_{M})$ is given by equation (\ref{Bateman_Horn_Const}).
\end{Conj}
We now apply this modified conjecture to two cases involving non-monic polynomials.

\subsection{Sophie Germain Primes} Application of the Bateman-Horn Conjecture to Sophie Germain primes was briefly discussed in \cite{BH02}, but the authors made no mention of its accuracy. A prime $p$ is a Sophie Germain prime if $2p+1$ is also a prime. Therefore we have two polynomials

\begin{align*}
f_{1}(n)&=n\\
f_{2}(n)&=2n+1
\end{align*}

For any prime $p>2$, the congruence equation $f(n)=n(2n+1)\equiv 0 \text{ (mod } p)$ has two solutions, corresponding to $n\equiv 0 \text{ (mod } p)$ and $n\equiv \frac{p-1}{2} \text{ (mod } p)$. For $p=2$, $f(n)=n(2n+1)\equiv 0\text{ (mod } p)$ has only one solution $n\equiv 0 \text{ (mod } p)$. Therefore the Bateman-Horn constant $C$ in equation (\ref{Bateman_Horn_Const}) is given by

\begin{align*}
C=2\prod_{p>2}\frac{p(p-2)}{(p-1)^{2}}=2C_{2}
\end{align*}
where $C_{2}=0.66016181584686957393$ \cite{Twinprime1}  is the Hardy-Littlewood twin-prime constant \cite{HL1}\cite{HL2}. Therefore equation (\ref{Modified_Conj_Equation}) becomes

\begin{equation}
\pi_{SG}(x)\sim 2C_{2}\int_{2}^{x}\frac{dt}{(\log t) \log(2t+1)}
\end{equation}

We generated about 29 million Sophie Germain primes using the freely available mathematical software SageMath v8.3 \cite{SageMath}. The largest Sophie Germain prime in the set is of the order $1.1\times 10^{10}$. The program uses the prime set generated by the function Primes(), and for every prime $p$ in the prime set, it tests whether $2p+1$ is also in the prime set. Table (\ref{Table-SophieGermain}) lists our results. The column with the heading $\pi_{SG}(x)$ corresponds to the actual counting from our empirical computations, while the next two columns correspond to the modified Bateman-Horn conjecture and the originally proposed conjecture, respectively. As can be clearly seen, the modified conjecture gives remarkably close results compared with the actual counting, while the original version of the conjecture produces rather significant deviations. We point out that  the modified version of the conjecture gives remarkably acurate estimates for even very modest values of $x$. That is, equation (\ref{Modified_Conj_Equation}) is not just an asymptotic expression, it seems to give very close estimates for almost all values of $x$. 

\begin{table}
\caption{$\pi_{SG}(x)$: Actual versus Conjectures}
\label{Table-SophieGermain}
\begin{tabular}{|c|c|c|c|}
\hline
$x$ & $\pi_{SG}(x)$ & $2C_{2}\int_{2}^{x}\frac{dt}{(\log t)\log (2t+1)}$ & $2C_{2} \int_{2}^{x}\frac{dt}{(\log t)^{2}}$\\
\hline\hline
$10^{2}$ & 10 & 10 & 14\\
\hline
$10^{3}$ & 37 & 39 & 46\\
\hline
$10^{4}$ & 190  & 195 & 214 \\
\hline
$10^{5}$ & 1171 & 1166 & 1249\\
\hline
$10^{6}$ & 7746 & 7811 & 8248 \\
\hline
$10^{7}$ & 56032 & 56128 & 58754 \\
\hline
$10^{8}$ & 423140 & 423294 & 440368\\
\hline
$10^{9}$& 3308859 & 3307888 & 3425308\\
\hline
$10^{10}$ & 26569515 & 26568824 & 27411417 \\
\hline
\end{tabular}
\end{table}

\subsection{$f(n)=6n^{2}+1$} 
We now present empirical data on the polynomial $f(n)=6n^{2}+1$. Denote the prime counting function by  $\pi_{6}(x)$. We need to calculate $\omega(p)$, the number of solutions to the congruence equation $6n^{2}+1\equiv 0 \text{ (mod }  p)$. For primes $p=2,3$, we have $\omega(2)=\omega(3)=0$. For any prime $p>3$, note that this polynomial has discriminant $D=-24$, therefore $\omega(p)=1+(-24\,|\,p)$ where $(-24\,|\,p)$ is the Legendre symbol. The Bateman-Horn constant in equation (\ref{Bateman_Horn_Const}) is then

\begin{equation}\label{HL-Const-6}
C=3\prod_{p> 3}\bigg(1-\frac{(-24\,|\,p)}{p-1}\bigg)
\end{equation}
which is essentially the expression given by the Hardy-Littlewood Conjecture F \cite{HL1}\cite{HL2}. This infinite product is known to converge extremely slowly (and not absolutely). One technique\footnote{This technique was suggested by K. Conrad on a Mathoverflow discussion forum: https://mathoverflow.net/questions/31150/calculating-the-infinite-product-from-the-hardy-littlewood-conjecture-f} to accelerate its convergence is to multiply the following identity on both sides of equation (\ref{HL-Const-6}):

\begin{equation}
L(1,(-24\,|\,\cdot))=\prod_{p> 3}\bigg(1-\frac{(-24\,|\,p)}{p}\bigg)^{-1}
\end{equation}
We then have

\begin{equation}\label{HL-Const-6a}
C=\frac{3}{L(1,(-24\,|\,\cdot))}\prod_{p> 3}\Bigg(\frac{{1-\frac{(-24\,|\, p)}{p-1}}}{1-\frac{(-24\,|\, p)}{p}}\Bigg)
\end{equation}
The infinite product in equation (\ref{HL-Const-6a}) converges as $p^{-2}$. We note that the character $(-24\,|\,\cdot)$ is primitive with a Conrey notation $\chi_{24}(5,\cdot)$ in the online $L$-function database \cite{LMFDB}, and the value of its corresponding Dirichlet $L$-function at $s=1$ has a closed form, $L(1,(-24\,|\,\cdot))=\pi/\sqrt{6}$. Truncating the infinite product in equation (\ref{HL-Const-6a}) up to primes $\sim 10^6$, we get $C=2.139124879$.

We generated 148 million primes of the form $6n^{2}+1$ using the mathematical software SageMath v8.3. The primality test is performed using the built-in function is$\_$prime(). The largest prime thus generated has the order of magnitude of $5.4\times 10^{19}$, and the corresponding value for $n\approx 7.3\times 10^{9}$. Table (\ref{Table-6N2+1}) lists the value of $\pi_{6}(x)$ and compares it with the estimates from the modified version of the conjecture as well as the original version of the conjecture. Again one can clearly see the remarkable accuracy of the modified conjecture even for modest values of $x$, while the original version of the conjecture shows considerable deviations, an indication that the value of $x\sim 10^{9}$ is too small for the asymptotic approximation to be accurate.

\begin{table}
\caption{$\pi_{6}(x)$: Actual versus Conjectures}
\label{Table-6N2+1}
\begin{tabular}{|c|c|c|c|}
\hline
$x$ & $\pi_{6}(x)$ & $C\int_{1}^{x}\frac{dt}{\log (6t^{2}+1)}$ & $\frac{C}{2} \int_{1}^{x}\frac{dt}{\log t}$\\
\hline\hline
$10^{2}$ & 27 & 25 & 31\\
\hline
$10^{3}$ & 155 & 162 & 189\\
\hline
$10^{4}$ & 1176 & 1195 & 1332\\
\hline
$10^{5}$ & 9445 & 9469 & 10299\\
\hline
$10^{6}$ & 78422 & 78514 & 84096 \\
\hline
$10^{7}$ & 671361 & 670963 & 711171\\
\hline
$10^{8}$ & 5859476 & 5859288 & 6163042\\
\hline
$10^{9}$ & 52007341 & 52009622 & 54386431\\
\hline
\end{tabular}
\end{table}
 
\section{Conclusion} Empirical computations suggest that the modified Bateman-Horn conjecture, equation (\ref{Modified_Conj_Equation}), to be quite accurate over a much wider range of prime values, which is strikingly similar to the fact that $\text{Li}(x)$ is a much better approximation of $\pi(x)$ than the asymptotic expression $\frac{x}{\log x}$ in the classic Prime Number Theorem. Further work is till needed, such as an estimate of error bounds, as well as possible adaptation to other problems such as the Mersenne primes.

\bibliographystyle{amsplain}

\end{document}